\documentclass{amsart}
\usepackage{hyperref}
 \newtheorem{theorem}{Theorem}[section]

 \newtheorem{conjecture}[theorem]{Conjecture}
 \theoremstyle{definition}
 
 \theoremstyle{remark}
 \newtheorem{remark}[theorem]{Remark}
 
 \numberwithin{equation}{section}

\newcommand{\C}{\mathbb{C}}
\newcommand{\R}{\mathbb{R}}
\newcommand{\W}{\mathbb{W}}

\newcommand{\eps}{\varepsilon}
\begin{document}

\title[Boundedness of the Cauchy Singular Integral Operator]
 {Remark on the Boundedness of the Cauchy Singular Integral Operator
 on Variable Lebesgue Spaces
 with Radial Oscillating Weights}
\author[Karlovich]{Alexei Yu. Karlovich}
\address{
Departamento de Matem\'atica,
Faculdade de Ci\^encias e Tecnologia,
Universidade Nova de Lisboa,
Quinta da Torre,
2829--516 Caparica,
Portugal}
\email{oyk@fct.unl.pt}
\thanks{The author is partially supported by F.C.T. (Portugal)
grant FCT/FEDER/POCTI/MAT /59972/2004.}


\date{\today}
\dedicatory{To Professor Kokilashvili on his seventieth birthday}
\begin{abstract}
Recently V. Kokilashvili, N. Samko, and S. Samko have proved a sufficient
condition for the boundedness of the Cauchy singular integral operator
on variable Lebesgue spaces with radial oscillating weights over Carleson
curves. This condition is formulated in terms of Matuszewska-Orlicz indices
of weights. We prove a partial converse of their result.
\end{abstract}

\maketitle
\section{Introduction and main result}
Let $\Gamma$ be a rectifiable curve in the complex plane. We equip $\Gamma$
with Lebesgue length measure $|d\tau|$. We say that a curve $\Gamma$ is simple
if it does not have self-intersections. In other words, $\Gamma$ is said to be simple
if it is homeomorphic
either to a line segment or to to a circle. In the latter situation we will say
that $\Gamma$ is a Jordan curve. The \textit{Cauchy singular integral} of
$f\in L^1(\Gamma)$ is defined by
\[
(Sf)(t):=\frac{1}{\pi i}\int_\Gamma\frac{f(\tau)}{\tau-t}d\tau
\quad(t\in\Gamma).
\]
This integral is understood in the principal value sense, that is,
\[
\int_\Gamma\frac{f(\tau)}{\tau-t}d\tau:=
\lim_{R\to 0}\int_{\Gamma\setminus\Gamma(t,R)}
\frac{f(\tau)}{\tau-t}d\tau,
\]
where $\Gamma(t,R):=\{\tau\in\Gamma:|\tau-t|<R\}$ for $R>0$.
David \cite{David84} (see also \cite[Theorem~4.17]{BK97}) proved that the
Cauchy singular integral generates the bounded operator $S$ on the Lebesgue
space $L^p(\Gamma)$, $1<p<\infty$, if and only if $\Gamma$ is a
\textit{Carleson} (\textit{Ahlfors-David regular}) \textit{curve}, that is,
\[
\sup_{t\in\Gamma}\sup_{R>0}\frac{|\Gamma(t,R)|}{R}<\infty,
\]
where for any measurable set $\Omega\subset\Gamma$ the symbol
$|\Omega|$ denotes its measure. To have a better idea about Carleson curves,
consider the following example. Let $\alpha>0$ and
\[
\Gamma:=\{0\}\cup\big\{\tau\in\C:\ \tau=x+ix^\alpha\sin(1/x),\ 0< x\le 1\big\}.
\]
One can show (see \cite[Example~1.3]{BK97}) that $\Gamma$ is not rectifiable
for $0<\alpha\le 1$, $\Gamma$ is rectifiable but not Carleson for $1<\alpha<2$,
and $\Gamma$ is a Carleson curve for $\alpha\ge 2$.

A measurable function $w:\Gamma\to[0,\infty]$ is referred to as a \textit{weight
function} or simply a \textit{weight} if $0<w(\tau)<\infty$ for almost all
$\tau\in\Gamma$. Suppose $p:\Gamma\to[1,\infty]$ is a measurable a.e. finite function.
Denote by $L^{p(\cdot)}(\Gamma,w)$ the set of all measurable complex-valued
functions $f$ on $\Gamma$ such that
\[
\int_\Gamma |f(\tau)w(\tau)/\lambda|^{p(\tau)}|d\tau|<\infty
\]
for some $\lambda=\lambda(f)>0$. This set becomes a Banach space when equipped
with the Luxemburg-Nakano norm
\[
\|f\|_{p(\cdot),w}:=\inf\left\{\lambda>0:
\int_\Gamma |f(\tau)w(\tau)/\lambda|^{p(\tau)}|d\tau|\le 1\right\}.
\]
If $p$ is constant, then $L^{p(\cdot)}(\Gamma,w)$ is nothing else but
the weighted Lebesgue space. Therefore, it is natural to refer to
$L^{p(\cdot)}(\Gamma,w)$ as a \textit{weighted generalized Lebesgue space
with variable exponent} or simply as a \textit{weighted variable Lebesgue
space}. This is a special case of Musielak-Orlicz spaces \cite{Musielak83}
(see also \cite{KR91}). Nakano \cite{Nakano50} considered these spaces
(without weights) as examples of so-called modular spaces, and sometimes
the spaces $L^{p(\cdot)}(\Gamma,w)$ are referred to as weighted Nakano
spaces.

Following \cite[Section~2.3]{KSS07}, denote by $W$ the class of all
continuous functions $\varrho: [0,|\Gamma|]\to[0,\infty)$ such that
$\varrho(0)=0$, $\varrho(x)>0$ if $0<x\le |\Gamma|$, and $\varrho$ is
almost increasing, that is, there is a universal constant $C>0$ such
that $\varrho(x)\le C\varrho(y)$ whenever $x\le y$. Further, let $\W$ be
the set of all functions $\varrho:[0,|\Gamma|]\to[0,\infty]$ such that
$x^\alpha\varrho(x)\in W$ and $x^\beta/\varrho(x)\in W$ for some
$\alpha,\beta\in\R$. Clearly, the functions $\varrho(x)=x^\gamma$
belong to $\W$ for all $\gamma\in\R$. For $\varrho\in\W$, put
\[
\Phi_\varrho^0(x):=\limsup_{y\to 0}\frac{\varrho(xy)}{\varrho(y)},
\quad
x\in(0,\infty).
\]
Since $\varrho\in\W$, one can show that the limits
\[
m(\varrho):=\lim_{x\to 0}\frac{\log\Phi_\varrho^0(x)}{\log x},
\quad
M(\varrho):=\lim_{x\to\infty}\frac{\log\Phi_\varrho^0(x)}{\log x}
\]
exist and $-\infty<m(\varrho)\le M(\varrho)<+\infty$. These numbers were
defined by Matuszewska
and Orlicz \cite{MO60,MO65} (see also \cite{M85} and \cite[Chapter~11]{M89}).
We refer to $m(\varrho)$ (resp. $M(\varrho)$) as the \textit{lower}
(resp. \textit{upper}) \textit{Matuszewska-Orlicz index} of $\varrho$.
For $\varrho(x)=x^\gamma$ one has $m(\varrho)=M(\varrho)=\gamma$.
Examples of functions $\varrho\in\W$ with $m(\varrho)<M(\varrho)$ can be
found, for instance, in \cite{AK89}, \cite[p.~93]{M89}, \cite[Section~2]{NSamko06}.

Fix pairwise distinct points $t_1,\dots,t_n\in\Gamma$ and functions
$w_1,\dots,w_n\in\W$. Consider the following weight
\begin{equation}\label{eq:weight}
w(t):=\prod_{k=1}^n w_k(|t-t_k|),
\quad
t\in\Gamma.
\end{equation}
Each function $w_k(|t-t_k|)$ is a radial oscillating weight.
The weight \eqref{eq:weight} is a continuous function on $\Gamma\setminus\{t_1,\dots,t_n\}$.
This is a
natural generalization of so-called Khvedelidze weights
$w(t)=\prod_{k=1}^n|t-t_k|^{\lambda_k}$, where $\lambda_k\in\R$
(see, e.g., \cite[Section~2.2]{BK97}, \cite{K56}, \cite{KPS06}).
Recently V. Kokilashvili, N. Samko, and S. Samko have proved the following
(see \cite[Theorem 4.3]{KSS07} and also \cite{KSS07-JFSA} for similar
results for maximal functions).
\begin{theorem}[V. Kokilashvili, S. Samko, N. Samko {\cite[Theorem 4.3]{KSS07}}]
\label{th:KSS}
Suppose $\Gamma$ is a simple Carleson curve and $p:\Gamma\to(1,\infty)$
is a continuous function satisfying
\begin{equation}\label{eq:Dini-Lipschitz}
|p(\tau)-p(t)|\le -A_\Gamma/\log|\tau-t|
\quad\mbox{whenever}\quad
|\tau-t|\le1/2,
\end{equation}
where $A_\Gamma$ is a positive constant depending only on $\Gamma$. Let
$w_1,\dots,w_n\in\W$ and the weight $w$ be given by \eqref{eq:weight}.
If
\begin{equation}\label{eq:KSS-condition}
0<1/p(t_k)+m(w_k),
\quad
1/p(t_k)+M(w_k)<1
\quad\mbox{for all}\quad
k\in\{1,\dots,n\},
\end{equation}
then the Cauchy singular integral operator $S$ is bounded on $L^{p(\cdot)}(\Gamma,w)$.
\end{theorem}
For the weight $w(t)=\prod_{k=1}^n|t-t_k|^{\lambda_k}$, \eqref{eq:KSS-condition}
reads as $0<1/p(t_k)+\lambda_k<1$ for all $k\in\{1,\dots,n\}$. This condition is
also necessary for the boundedness of $S$ on the variable Lebesgue space
$L^{p(\cdot)}(\Gamma,w)$ with the Khvedelidze weight $w$ (see \cite{KPS06}).

The author have proved in \cite{Karlovich-preprint}
that for Jordan curves condition \eqref{eq:KSS-condition} is necessary
for the boundedness of the operator $S$.
\begin{theorem}[{\cite[Corollary~4.3]{Karlovich-preprint}}]
Suppose $\Gamma$ is a rectifiable Jordan curve and $p:\Gamma\to(1,\infty)$
is a continuous function satisfying \eqref{eq:Dini-Lipschitz}.
Let $w_1,\dots,w_n\in\W$ and the weight $w$ be given by \eqref{eq:weight}.
If the Cauchy singular integral operator $S$ is bounded on $L^{p(\cdot)}(\Gamma,w)$,
then $\Gamma$ is a Carleson curve and \eqref{eq:KSS-condition} is fulfilled.
\end{theorem}
The proof of this result given in \cite{Karlovich-preprint} essentially uses
that $\Gamma$ is closed. In this paper we embark on the situation
of non-closed curves. Our main result is a partial converse of Theorem~\ref{th:KSS}.
It follows from our results \cite{Karlovich03,Karlovich-preprint} based
on further development of ideas from \cite[Chap.~1--3]{BK97}.
\begin{theorem}[Main result]
\label{th:main}
Let $\Gamma$ be a rectifiable curve homeomorphic to a line segment and $p:\Gamma\to(1,\infty)$
be a continuous function satisfying \eqref{eq:Dini-Lipschitz}.
Suppose $w_1,\dots,w_n\in\W$ and the weight $w$ is given by \eqref{eq:weight}.
If the Cauchy singular integral operator $S$ is bounded on $L^{p(\cdot)}(\Gamma,w)$,
then $\Gamma$ is a Carleson curve and
\[
0\le 1/p(t_k)+m(w_k),
\quad
1/p(t_k)+M(w_k)\le 1
\quad\mbox{for all}\quad
k\in\{1,\dots,n\}.
\]
Moreover, if there exists an $\eps_0>0$ such that the Cauchy singular
integral operator $S$ is bounded on $L^{p(\cdot)}(\Gamma,w^{1+\eps})$
for all $\eps\in(-\eps_0,\eps_0)$, then
\[
0< 1/p(t_k)+m(w_k),
\quad
1/p(t_k)+M(w_k)< 1
\quad\mbox{for all}\quad
k\in\{1,\dots,n\}.
\]
\end{theorem}
For standard Lebesgue spaces, the boundedness of the operator $S$ on
$L^p(\Gamma,w)$, $1<p<\infty$, implies that $S$ is also bounded on
$L^p(\Gamma,w^{1+\eps})$ for all $\eps$ in a sufficiently small neighborhood
of zero (see \cite[Theorems~2.31 and 4.15]{BK97}). Hence if $1<p<\infty$,
$\Gamma$ is a simple Carleson curve, $w_1,\dots,w_n\in\W$, and
the weight $w$ is given by \eqref{eq:weight}, then $S$ is bounded
on the standard Lebesgue space $L^p(\Gamma,w)$, $1<p<\infty$, if and only
if
\[
0< 1/p+m(w_k),
\quad
1/p+M(w_k)< 1
\quad\mbox{for all}\quad
k\in\{1,\dots,n\}.
\]
We believe that all weighted variable Lebesgue spaces have this stability property.
\begin{conjecture}\label{co:stability}
Let $\Gamma$ be a simple rectifiable curve, $p:\Gamma\to[1,\infty]$ be a measurable
a.e. finite function, and $w:\Gamma\to[0,\infty]$ be a weight such that the Cauchy
singular integral operator $S$ is bounded on $L^{p(\cdot)}(\Gamma,w)$. Then
there is a number $\eps_0>0$ such that $S$ is bounded on $L^{p(\cdot)}(\Gamma,w^{1+\eps})$
for all $\eps\in(-\eps_0,\eps_0)$.
\end{conjecture}

If this conjecture would be true, we were able to prove the complete converse
of Theorem~\ref{th:KSS} for non-closed curves, too.
\section{Proof}
In this section we formulate several results from
\cite{BK97,Karlovich03,Karlovich-preprint}
and show that Theorem~\ref{th:main} easily follows from them.
\subsection{Muckenhoupt type condition}
Suppose $\Gamma$ is a simple rectifiable curve and $p:\Gamma\to(1,\infty)$ is
a continuous function. Since $\Gamma$ is compact, one has
\[
1<\min_{\tau\in\Gamma} p(\tau),
\quad
\max_{\tau\in\Gamma} p(\tau)<\infty
\]
and the conjugate exponent
\[
q(\tau):=p(\tau)/(p(\tau)-1)\quad(\tau\in\Gamma)
\]
is well defined and also bounded and bounded away from zero.
We say that a weight $w:\Gamma\to[0,\infty]$ belongs
to $A_{p(\cdot)}(\Gamma)$ if
\[
\sup_{t\in\Gamma}\sup_{R>0}\frac{1}{R}
\|w\chi_{\Gamma(t,R)}\|_{p(\cdot)}
\|w^{-1}\chi_{\Gamma(t,R)}\|_{q(\cdot)}<\infty.
\]
If $p=const\in(1,\infty)$, then
this class coincides with the well known Muckenhoupt class. From the H\"older
inequality for $L^{p(\cdot)}(\Gamma)$ (see e.g. \cite[Theorems 13.12 and 13.13]{Musielak83}
for Muslielak-Orlicz spaces over arbitrary measure spaces and also
\cite[Theorem~2.1]{KR91} for variable Lebesgue spaces over domains in $\R^n$) it
follows that if $w\in A_{p(\cdot)}(\Gamma)$, then $\Gamma$ is a Carleson curve.

Since $L^{p(\cdot)}(\Gamma,w)$ is a Banach function space in the sense of
\cite[Definition~1.1]{BS88}, the next result follows from
\cite[Theorem~6.1]{Karlovich03} (stated in \cite{Karlovich03}
for Jordan curves, however its proof remains the same for curves
homeomorphic to line segments, see also \cite[Theorem~3.2]{Karlovich98}).
\begin{theorem}\label{th:Ap}
Let $\Gamma$ be a simple rectifiable curve and let $p:\Gamma\to(1,\infty)$
be a continuous function. If $w:\Gamma\to[0,\infty]$ is an arbitrary weight
such that the operator $S$ is bounded on $L^{p(\cdot)}(\Gamma,w)$, then
$w\in A_{p(\cdot)}(\Gamma)$.
\end{theorem}
If $p=const\in(1,\infty)$, then $w\in A_p(\Gamma)$ is also sufficient
for the boundedness of $S$ on the weighted Lebesgue space $L^p(\Gamma,w)$
(see e.g. \cite[Theorem~4.15]{BK97}).
\subsection{Submultiplicative functions}
Following \cite[Section~1.4]{BK97}, we say a function
$\Phi:(0,\infty)\to(0,\infty]$ is \textit{regular} if it is bounded in an
open neighborhood of $1$. A function $\Phi:(0,\infty)\to(0,\infty]$ is said
to be \textit{submultiplicative} if
\[
\Phi(xy)\le\Phi(x)\Phi(y)
\quad\mbox{for all}\quad x,y\in(0,\infty).
\]
It is easy to show that if $\Phi$ is regular and submultiplicative, then
$\Phi$ is bounded away from zero in some open neighborhood of $1$.
Moreover, in this case $\Phi(x)$ is finite for all $x\in(0,\infty)$.
Given a regular and submultiplicative function $\Phi:(0,\infty)\to(0,\infty)$,
one defines
\[
\alpha(\Phi):=\sup_{x\in(0,1)}\frac{\log\Phi(x)}{\log x},
\quad
\beta(\Phi):=\inf_{x\in(1,\infty)}\frac{\log\Phi(x)}{\log x}.
\]
Clearly, $-\infty<\alpha(\Phi)$ and $\beta(\Phi)<\infty$.
\begin{theorem}[{see \cite[Theorem~1.13]{BK97} or \cite[Chap.~2, Theorem~1.3]{KPS82}}]
\label{th:submult}
If a function $\Phi:(0,\infty)\to(0,\infty)$ is regular and submultiplicative, then
\[
\alpha(\Phi)=\lim_{x\to 0}\frac{\log\Phi(x)}{\log x},
\quad
\beta(\Phi)=\lim_{x\to\infty}\frac{\log\Phi(x)}{\log x}
\]
and $-\infty<\alpha(\Phi)\le\beta(\Phi)<+\infty$.
\end{theorem}
The quantities $\alpha(\Phi)$ and $\beta(\Phi)$ are called the \textit{lower}
and  \textit{upper indices of the regular and submultiplicative function}
$\Phi$, respectively.
\subsection{Indices of powerlikeness}
Fix $t\in\Gamma$ and put $d_t:=\max\limits_{\tau\in\Gamma}|\tau-t|$.
Suppose $w:\Gamma\to[0,\infty]$ is a weight such that $\log w\in L^1(\Gamma(t,R))$ for
every $R\in(0,d_t]$. Put
\[
H_{w,t}(R_1,R_2):=
\frac{\displaystyle\exp\left(\frac{1}{|\Gamma(t,R_1)|}
\int_{\Gamma(t,R_1)}\log w(\tau)|d\tau|\right)}
{\displaystyle\exp\left(\frac{1}{|\Gamma(t,R_2)|}
\int_{\Gamma(t,R_2)}\log w(\tau)|d\tau|\right)},
\quad R_1,R_2\in(0,d_t].
\]
Consider the function
\[
(V_t^0w)(x) := \limsup_{R\to 0} H_{w,t}(xR,R),
\quad x\in(0,\infty).
\]
Combining Lemmas~4.8--4.9 and Theorem~5.9 of \cite{Karlovich03} with
Theorem~3.4, Lemma~3.5 of \cite{BK97}, we arrive at the following.
\begin{theorem}\label{th:indices-inequality}
Let $\Gamma$ be a simple rectifiable curve, $p:\Gamma\to(1,\infty)$
be a continuous function satisfying \eqref{eq:Dini-Lipschitz}, and
$w:\Gamma\to[0,\infty]$ be a weight such that $w\in A_{p(\cdot)}(\Gamma)$.
Then, for every $t\in\Gamma$, the function $V_t^0 w$ is regular and
submultiplicative and
\[
0\le 1/p(t)+\alpha(V_t^0w),
\quad
1/p(t)+\beta(V_t^0w)\le 1.
\]
\end{theorem}
The numbers $\alpha(V_t^0w)$ and $\beta(V_t^0w)$ are called the
\textit{lower} and \textit{upper indices of powerlikeness} of $w$ at
$t\in\Gamma$, respectively (see \cite[Chap.~3]{BK97}). This terminology
can be explained by the simple fact that for the power weight
$w(\tau):=|\tau-t|^\lambda$ its indices of powerlikeness coincide and are equal
to $\lambda$.
\subsection{Matuszewska-Orlicz indices as indices of powerlikeness}
If $\varrho\in\W$, then $\Phi_\varrho^0$ is a regular and submultiplicative
function and its indices are nothing else but the Matuszewska-Orlicz indices
$m(\varrho)$ and $M(\varrho)$. The next result shows that for radial
oscillating weights indices of powerlikeness and Matuszewska-Orlicz indices
coincide.
\begin{theorem}[{see \cite[Theorem 2.8]{Karlovich-preprint}}]
\label{th:radial-MO}
Suppose $\Gamma$ is a simple Carleson curve. If $w_1,\dots,w_n\in\W$
and $w(\tau)=\prod_{k=1}^n w_k(|\tau-t_k|)$, then for every $t\in\Gamma$
the function $V_t^0w$ is regular and submultiplicative and
\[
\begin{array}{llll}
\alpha(V_{t_k}^0w)=m(w_k), & \beta(V_{t_k}^0w)=M(w_k) &\mbox{for}  &k\in\{1,\dots,n\},
\\[2mm]
\alpha(V_t^0w)=0, & \beta(V_t^0w)=0 &  \mbox{for}   &t\in\Gamma\setminus\{t_1,\dots,t_n\}.
\end{array}
\]
\end{theorem}
Note that in \cite{Karlovich-preprint}, Theorem~\ref{th:radial-MO} is proved
for Jordan curves. But the proof does not use the assumption that $\Gamma$
is closed. It works also for non-closed curves considered in this paper.
\subsection{Proof of Theorem~\ref{th:main}}
\begin{proof}
Suppose $S$ is bounded on $L^{p(\cdot)}(\Gamma,w)$. From Theorem~\ref{th:Ap}
it follows that $w\in A_{p(\cdot)}(\Gamma)$. By H\"older's inequality
this implies that $\Gamma$ is a Carleson curve. Fix an arbitrary $t\in\Gamma$.
Then, in view of
Theorems~\ref{th:submult} and \ref{th:indices-inequality} the function
$V_t^0w$ is regular and submultiplicative, so its indices are well defined and
satisfy $0\le 1/p(t)+\alpha(V_t^0w)$ and $1/p(t)+\beta(V_tw)\le 1$.
From these inequalities and Theorem~\ref{th:radial-MO} it follows that
\begin{equation}\label{eq:basic}
0\le 1/p(t_k)+m(w_k),
\quad
1/p(t_k)+M(w_k)\le 1
\end{equation}
for all $k\in\{1,\dots,n\}$.

If $S$ is bounded on all spaces $L^{p(\cdot)}(\Gamma,w^{1+\eps})$ for all
$\eps$ in a neighborhood of zero, then as before
\[
0\le 1/p(t_k)+m(w_k^{1+\eps}),\quad 1/p(t_k)+M(w_k^{1+\eps})\le 1
\]
for every $k\in\{1,\dots,n\}$. It is easy to see that $m(w_k^{1+\eps})=(1+\eps)m(w_k)$
and $M(w_k^{1+\eps})=(1+\eps)M(w_k)$. Therefore
\[
0\le 1/p(t_k) +(1+\eps)m(w_k),
\quad
1/p(t_k)+(1+\eps)M(w_k)\le 1
\]
for all $\eps$ in a neighborhood of zero and for all $k\in\{1,\dots,n\}$.
These inequalities immediately imply that $0<1/p(t_k)+m(w_k)$
and $1/p(t_k)+M(w_k)<1$ for all $k$.
\end{proof}
\begin{remark}
The presented proof involves the notion of indices of powerlikeness,
which were invented to treat general Muckenhoupt weights (see \cite{BK97}).
Weights considered in the present paper are continuous except for a finite
number of points. So, it would be rather interesting to find a direct proof
of the fact that $w\in A_{p(\cdot)}(\Gamma)$ implies \eqref{eq:basic}, which
does not involve the indices of powerlikeness $\alpha(V_t^0w)$ and $\beta(V_t^0w)$.
\end{remark}
\subsection{Final remarks}
In connection with Conjecture~\ref{co:stability}, we would like to note that
for standard Lebesgue spaces $L^p(\Gamma,w)$ there are two different proofs of
the stability of the boundedness of $S$ on $L^p(\Gamma,w^{1+\eps})$ for small
$\eps$. Simonenko's proof \cite{Simonenko83} is based on the stability of the
Fredholm property of some singular integral operators related to the Riemann
boundary value problem. Another proof is based on the self-improving property
of Muckenhoupt weights (see e.g. \cite[Theorem~2.31]{BK97}). One may ask
whether does $w\in A_{p(\cdot)}(\Gamma)$ imply $w^{1+\eps}\in A_{p(\cdot)}(\Gamma)$
for all $\eps\in(-\eps_0,\eps_0)$ with some fixed $\eps>0$? The positive answer
would give a proof of the complete converse of \eqref{eq:KSS-condition}.
The author does not know any stability result for the boundedness of $S$
or a self-improving property for $w\in A_{p(\cdot)}(\Gamma)$.

After this paper had been submitted, P.~H\"ast\"o and L.~Diening \cite{HD08}
have found a necessary and sufficient condition for the
boundedness of the classical Hardy-Littlewood maximal function on
weighted variable Lebesgue spaces in the setting of $\R^n$. Note
that they write a weight as a measure (outside of
$|\cdot|^{p(\tau)}$). Their condition is another generalization of
the classical Muckenhoupt condition.
In the setting of Carleson curves (and
the weight written inside of $|\cdot|^{p(\tau)}$), the
H\"ast\"o-Diening condition takes the form
\begin{equation}\label{eq:Hasto-Diening}
\sup_{t\in\Gamma}\sup_{R>0}\left(\frac{1}{R^{p_{\Gamma(t,R)}}}
\int_{\Gamma(t,R)}w(\tau)^{p(\tau)}|d\tau|\right)
\|w(\cdot)^{-p(\cdot)}\chi_{\Gamma(t,R)}(\cdot)\|_{q(\cdot)/p(\cdot)}<\infty,
\end{equation}
where
\[
p_{\Gamma(t,R)}:=\left(\frac{1}{|\Gamma(t,R)|}\int_{\Gamma(t,R)}\frac{1}{p(\tau)}\,|d\tau|\right)^{-1}.
\]
Let $HD_{p(\cdot)}(\Gamma)$ denote the class of weights
$w:\Gamma\to[0,\infty]$ satisfying \eqref{eq:Hasto-Diening}.
Following the arguments contained in \cite[Remark~3.10]{HD08},
one can show that
\[
A_{L^{p(\cdot)}}(\Gamma)\supset HD_{p(\cdot)}(\Gamma)
\]
whenever $p:\Gamma\to(1,\infty)$ satisfies the Dini-Lipschitz condition
\eqref{eq:Dini-Lipschitz}. We conjecture that the
H\"ast\"o-Diening characterization remains true also for the operator $S$ in the
setting of Carleson curves.
\begin{conjecture}\label{co:Hasto-Diening}
Let $\Gamma$ be a simple Carleson curve, $w:\Gamma\to[0,\infty]$ be a
weight, and $p:\Gamma\to(1,\infty)$ be a continuous function
satisfying the Dini-Lipschitz condition \eqref{eq:Dini-Lipschitz}.
The operator $S$ is bounded on $L^{p(\cdot)}(\Gamma,w)$ if and
only if $w\in HD_{p(\cdot)}(\Gamma)$.
\end{conjecture}
\subsection*{Acknowledgement.} The author would like to thank the anonymous
referee for several useful remarks.


\begin{thebibliography}{99}
\bibitem{AK89}
V. D. Aslanov and Yu. I. Karlovich,
\textit{One-sided invertibility of functional operators in reflexive Orlicz spaces.}
Akad. Nauk Azerbaidzhan. SSR Dokl. \textbf{45} (1989), no. 11--12, 3--7 (in Russian).

\bibitem{BS88}
C. Bennett and R. Sharpley,
\textit{Interpolation of Operators.}
Academic Press, Boston, 1988.

\bibitem{BK97}
A. B\"ottcher and Yu. I. Karlovich,
\textit{Carleson Curves, Muckenhoupt Weights, and Toeplitz Operators.}
Birkh\"auser, Basel, 1997.

\bibitem{David84}
G. David,
\textit{Oper\'ateurs int\'egraux singuliers sur certaines courbes du plan
complexe.}
Ann. Sci. \'Ecole Norm. Super. \textbf{17} (1984), 157--189.

\bibitem{HD08}
P.~H\"ast\"o and L.~Diening, Muckenhoupt weights in variable
exponent spaces, (December 15, 2008), preprint is available at
{http://www.helsinki.fi/$\widetilde{\hspace{2mm}}$hasto/pp/p$75\underline{\hspace{2mm}}$submit.pdf}.

\bibitem{Karlovich03}
A. Yu. Karlovich,
\textit{Fredholmness of singular integral operators with piecewise continuous
coefficients on weighted Banach function spaces.}
J. Integr. Equat. Appl. \textbf{15} (2003), 263--320.

\bibitem{Karlovich98}
A. Yu. Karlovich,
\textit{Singular integral operators with piecewise continuous coefficients
in reflexive rearrangement-invariant spaces.}
Integral Equatations Operator Theory \textbf{32} (1998), 436--481.

\bibitem{Karlovich-preprint}
A. Yu. Karlovich,
\textit{Singular integral operators on variable Lebesgue spaces
with radial oscillating weights.}
Proceedings of IWOTA 2007, Operator Theory: Advances and Applications, Birk¨h\"auser, to appear.
Preprint is available at {\tt http://arxiv.org/abs/0708.0778}.

\bibitem{K56}
B. V. Khvedelidze,
\textit{Linear discontinuous boundary problems in the theory of functions,
singular integral equations and some of their applications.}
Trudy Tbiliss. Mat. Inst. Razmadze \textbf{23} (1956), 3--158 (in Russian).

\bibitem{KPS06}
V. Kokilashvili, V. Paatashvili, and S. Samko,
\textit{Boundedness in Lebesgue spaces with variable exponent of the Cauchy
singular operator on Carleson curves.}
In: ``Modern Operator Theory and Applications. The Igor Borisovich Simonenko
Anniversary Volume".
Operator Theory: Advances and Applications \textbf{170} (2006), 167--186.

\bibitem{KSS07-JFSA}
V. Kokilashvili, N. Samko, and S. Samko,
\textit{The maximal operator in weighted variable spaces $L^{p(\cdot)}$.}
J. Func. Spaces Appl. \textbf{5} (2007), 299--317.

\bibitem{KSS07}
V. Kokilashvili, N. Samko, and S. Samko,
\textit{Singular operators in variable spaces $L^{p(\cdot)}(\Omega,\rho)$
with oscillating weights.}
Math. Nachr. \textbf{280} (2007), 1145--1156.

\bibitem{KR91}
O. Kov\'a{\v c}ik and J.~R\'akosn{\'\i}k,
\textit{On spaces $L\sp {p(x)}$ and $W\sp {k,p(x)}$.}
Czechoslovak Math. J. \textbf{41(116)} (1991), 592--618.

\bibitem{KPS82}
S. G. Krein, Ju. I. Petunin, and E. M. Semenov,
\textit{Interpolation of Linear Operators.}
AMS Translations of Mathematical Monographs \textbf{54},
Providence, RI, 1982.

\bibitem{M85}
L. Maligranda,
\textit{Indices and interpolation.}
Dissert. Math. \textbf{234} (1985), 1--49.

\bibitem{M89}
L. Maligranda,
\textit{Orlicz Spaces and Interpolation.}
Sem. Math. 5, Dep. Mat.,
Universidade Estadual de Campinas, Campinas SP, Brazil, 1989.

\bibitem{MO60}
W. Matuszewska and W. Orlicz,
\textit{On certain properties of $\varphi $-functions.}
Bull. Acad. Polon. Sci. Ser. Sci. Math. Astronom. Phys. \textbf{8} (1960), 439--443.
Reprinted in: W. Orlicz, \textit{Collected Papers}, PWN, Warszawa, 1988, 1112--1116.

\bibitem{MO65}
W. Matuszewska and W. Orlicz,
\textit{On some classes of functions with regard to their orders of growth.}
Studia Math. \textbf{26} (1965), 11--24.
Reprinted in: W. Orlicz, \textit{Collected Papers}, PWN, Warszawa, 1988, 1217--1230.

\bibitem{Musielak83}
J. Musielak,
\textit{Orlicz Spaces and Modular Spaces.}
Lecture Notes in Mathematics \textbf{1034}.
Springer-Verlag, Berlin, 1983.

\bibitem{Nakano50}
H. Nakano,
\textit{Modulared Semi-Ordered Linear Spaces.}
Maruzen Co., Ltd., Tokyo, 1950.

\bibitem{NSamko06}
N. Samko,
\textit{Singular integral operators in weighted spaces of continuous
functions with oscillating continuity moduli and oscillating weights.}
In: ``The Extended Field of Operator Theory".
Operator Theory: Advances and Applications \textbf{171} (2006), 323--347.

\bibitem{Simonenko83}
I. B. Simonenko,
\textit{Stability of the weight properties of functions with respect to a singular Cauchy integral},
Math Notes \textbf{33} (1983), 208--212.
\end{thebibliography}
\end{document}